\documentclass[12pt]{amsart}
\renewcommand{\bar}{\overline}
\usepackage{amsmath,amsthm,amscd,euscript}
\renewcommand{\bar}{\overline}
\def \r{\mathbb R}
\def \c{\mathbb C}
\def \q{\mathbb Q}
\def \z{\mathbb Z}

\DeclareMathOperator{\sign}{sign}
\DeclareMathOperator{\sspan}{Span}
\DeclareMathOperator{\il}{I\ell} 
\DeclareMathOperator{\isin}{lsin}

\newtheorem{theorem}{Theorem}[section]
\newtheorem{lemma}[theorem]{Lemma}

\newtheorem{proposition}[theorem]{Proposition}
\newtheorem{corollary}[theorem]{Corollary}
\theoremstyle{remark}
\newtheorem{remark}[theorem]{Remark}
\theoremstyle{definition}
\newtheorem{definition}[theorem]{Definition}
\newtheorem{example}[theorem]{Example}

\title{Multidimensional Gauss Reduction Theory for conjugacy classes of $SL(n,\z)$}
\author{Oleg~Karpenkov}
\date{02 May 2012}

\keywords{Gauss Reduction Theory, Klein-Voronoi continued
fractions, convex hulls, Hessenberg matrices, Markoff-Davenport
characteristic}

\email[Oleg Karpenkov]{karpenk@mccme.ru}
\address{TU Graz /Kopernikusgasse 24, A 8010 Graz, Austria/}

\begin{document}
\input epsf

\begin{abstract}
In this paper we describe the set of conjugacy classes in the
group $SL(n,\z)$. We expand geometric Gauss Reduction Theory that
solves the problem for $SL(2,\z)$ to the multidimensional case.
Further we find complete invariant of classes in terms of
multidimensional Klein-Voronoi continued fractions, where
$\varsigma$-reduce Hessenberg matrices play the role of reduced
matrices.
\end{abstract}

\maketitle
\tableofcontents

\section*{Introduction}\label{intro}

Two matrices $M_1$ and $M_2$ in $SL(n,\z)$ are {\it integer
conjugate} if there exists a matrix $X$ in $GL(n,\z)$ such that
$$
M_2=XM_1X^{-1}.
$$
 In this paper we study the following problem.

\vspace{2mm}

{\bf Problem.} Describe the set of integer conjugacy classes in
$SL(n,\z)$.

\vspace{2mm}

One of the mostly common strategies to solve this kind of problems
is to find complete invariants to distinguish the classes, and
further if possible to write normal form of conjugacy classes. For
instance, in the similar problem for $SL(n,F)$ for an
algebraically closed field $F$ one has Jordan Normal Forms as a
complete description of conjugacy classes. Jordan blocks form a
complete invariant in this case.

If the field is not algebraically closed, the description is much
more complicated via Jordan-Chevalley decomposition. In the study
of $SL(n,\z)$ we are faced with a group instead of a field. For
the general case it is only known the solution of the similarity
problem on verification whether two matrices are conjugate or not
(see in~\cite{Appelgate1982} and~\cite{Grunewald1980}). A complete
description of the set of integer conjugacy classes in $SL(2,\z)$
is given by Gauss Reduction Theory (see for instance
in~\cite{Katok2003} and \cite{Manin2002}). It turns out that it is
natural to consider several normal forms for an integer conjugacy
class but not necessarily only one.

Currently the main approach to the study of the above problem is
as follows: one should try to split $GL(n,\q)$-conjugacy classes
into $GL(n,\z)$ conjugacy classes. Then the problem is reduced to
certain problems related to orders of algebraic fields defined by
the roots of characteristic polynomial of the corresponding
matrices (like computing their class numbers, etc.). In this paper
we introduce an alternative geometric approach based on
generalization of Gauss Reduction Theorem. We will study questions
related to three-dimensional case in more details in our
forthcoming paper.

\vspace{2mm}

{\bf Description of the paper.} In current paper we present the
following main four results.

\vspace{2mm}

{\bf I. Matrix description of integer conjugacy classes.} We
consider Hessenberg matrices as a multidimensional analog of
reduced matrices in Gauss Reduction Theory. Hessenberg matrices
are matrices that vanish below the superdiagonal (for more
information see in~\cite{Stoer2002}). These matrices were
essentially used in the QR-algorithm for eigenvalue problem, but
they were never considered before in the frames of similarity
theory. We introduce a natural notion of Hessenberg complexity for
Hessenberg matrices, which is a nonnegative integer function, and
show that {\it each integer conjugacy class of irreducible
matrices has only finite number of Hessenberg matrices with
minimal complexity}. This result is a combination of
Theorem~\ref{toRedHess} and Theorem~\ref{finiteness}. We study all
related questions in Section~\ref{HMCC}.

\vspace{2mm}

{\bf II. Geometric complete invariant of integer conjugacy
classes.} In Section~\ref{KlVor} we introduce the complete
invariant of integer conjugacy classes of $GL(n,\z)$ matrices.
Recently in~\cite{Karpenkov2010} we showed a geometric explanation
of Gauss Reduction Theory in terms of geometric continued
fractions (see also briefly in Subsection~\ref{GRT}). In
Sections~\ref{KlVor} we extend this approach to the
multidimensional case. We propose a geometric description of
integer conjugacy classes in terms of multidimensional continued
fractions in the sense of Klein-Voronoi: {\it the periods of such
Klein-Voronoi continued fractions are complete invariants of
integer conjugacy classes} (Theorem~\ref{teorem_inv}). In addition
we study the group structure of the set of periods
(Theorem~\ref{PeriodStructure}).

\vspace{2mm}

{\bf III. Techniques to construct reduced matrices.} In
Section~\ref{HMCC} we introduce a techniques to construct
$\varsigma$-reduced matrices integer conjugate to a given one. It
is based on the following result (Theorem~\ref{cfr}): {\it any
$\varsigma$-reduced matrix is obtained from an integer vertex of
Klein-Voronoi continued fraction by applying to it the algorithm
of Subsection~\ref{toHess}}.

\vspace{2mm}

%

{\bf Acknowledgment.} The work is partially supported by FWF grant
M~1273-N18. The author is grateful to H.~W.~Lenstra and
E.~I.~Pav\-lov\-skaya for useful remarks.

\section{Hessenberg matrices and conjugacy classes}\label{HMCC} In this
section we study questions of reduction to $\varsigma$-reduced
matrices and investigate families of perfect Hessenberg matrices
in general. We start with necessary definitions and notation in
Subsection~\ref{notions}. In Subsection~\ref{Hess1} we prove that
any integer conjugacy class with irreducible characteristic
polynomial has at least one $\varsigma$-reduced matrix. Further in
Subsection~\ref{Hess2} we show that Hessenberg matrices are
defined by their Hessenberg types and characteristic polynomials
and deduce the finiteness of $\varsigma$-reduced matrices in each
integer conjugacy class. Finally in Subsection~\ref{Struc} we
study the structure of the set of perfect Hessenberg matrices. We
consider this set as a ''book'' that contains ''pages'' enumerated
by Hessenberg type. The matrices of the same page are
distinguished by characteristic polynomial, only matrices from
different pages can be integer conjugate.

\subsection{Notions and definition}\label{notions}

In this subsection we briefly introduce matrices that generalize
the reduced matrices in Gauss Reduction Theory for $SL(2,\z)$.

\vspace{2mm}

\subsubsection{Perfect Hessenberg matrices} A matrix $M$ of the
form
$$
\left(
\begin{array}{cccccc}
a_{1,1} &a_{1,2} &\cdots &a_{1,n-2}  &a_{1,n-1}  &a_{1,n}\\
a_{2,1} &a_{2,2} &\cdots &a_{2,n-2}  &a_{2,n-1}  &a_{2,n}\\
0       &a_{3,2} &\cdots &a_{3,n-2}  &a_{3,n-1}  &a_{3,n}\\
\vdots  &\vdots  &\ddots &\vdots     &\vdots     &\vdots \\
0       &0       &\cdots &a_{n-1,n-2}&a_{n-1,n-1}&a_{n-1,n}\\
0       &0       &\cdots &0          &a_{n,n-1}  &a_{n,n}\\
\end{array}
\right)
$$
is called an {\it $($upper$)$ Hessenberg} matrix. Such matrices
were first studied by K.~Hessenberg in~\cite{Hessenberg1942} and
further used in QR-algorithms (see~\cite{Demmel1997},
\cite{Trefethen1997}, and \cite{Ortega1963/1964}). We say that the
matrix $M$ is of {\it Hessenberg type}
$$
\langle
a_{1,1},a_{2,1}|a_{1,2},a_{2,2},a_{3,2}|\cdots|a_{1,n{-}1},
a_{2,n{-}1},\ldots,a_{n,n{-}1} \rangle.
$$

\begin{definition}
A Hessenberg matrix in $SL(n,\z)$ is said to be {\it perfect} if
for any pair of integers $(i,j)$ satisfying $1\le i<j{+}1\le n$
the following inequalities hold: $0\le a_{i,j}<a_{j+1,j}$.
\end{definition}

In other words all elements of all the columns except the last
column of a perfect Hessenberg matrix are nonnegative integers,
the maximal elements in these columns are the lowest nonzero ones
(i.e., $a_{j+1,j}$, $j=1,\ldots, n{-}1$).

\vspace{2mm}

\subsubsection{$\varsigma$-reduced Hessenberg matrices}  We mostly
study $SL(n,\z)$-matrices with irreducible characteristic
polynomials over rational numbers. Any such matrix is integer
conjugate to a perfect Hessenberg matrix with positive Hessenberg
complexity (see Theorem~\ref{toRedHess} below). Actually, there
are infinitely many perfect Hessenberg matrices integer conjugate
to a given one. So we give an additional notion of complexity to
reduce the number of such matrices.

\begin{definition}The integer number
$$
\prod\limits_{j=1}^{n{-}1}|a_{j{+}1,j}|^{n-j}
$$
is called the {\it Hessenberg complexity} of the matrix $M$ and
denoted by $\varsigma(M)$.
\end{definition}

Hessenberg complexity has the following geometric meaning. It is
equivalent to the volume of the parallelepiped spanned by $e_1$,
$M(e_1)$, $M^2(e_1), \ldots, M^{n-1}(e_1)$, where $e_1$ is the
first basis vector i.e.~$(1,0,\ldots, 0)$, we discuss this in more
details in Subsection~\ref{cf_2}.

An integer Hessenberg matrix has the unit Hessenberg complexity if
and only if $a_{2,1}=\cdots=a_{n,n-1}=1$, such matrices are called
{\it Frobenius} matrices. The elements of the last column of a
Frobenius matrix are the coefficients of the characteristic
polynomial multiplied alternatively by $\pm 1$.

\begin{example}
The following matrix
$$
\left(
\begin{array}{ccc}
1&2&2\\
2&3&4\\
0&5&-1\\
\end{array}
\right)
$$
is a perfect Hessenberg matrix of type $\langle 1,2|2,3,5
\rangle$. The Hessenberg complexity of this matrix is $2^2\cdot
5=20$.
\end{example}

\begin{definition}
We say that a perfect Hessenberg matrix $M$ is {\it
$\varsigma$-reduced} if its Hessenberg complexity is the least
possible. Otherwise we say that the matrix is {\it
$\varsigma$-nonreduced}.
\end{definition}

In Theorem~\ref{finiteness} we show that the number of
$\varsigma$-reduced matrices is finite in any integer conjugacy
class. Still sometimes there are several $\varsigma$-reduced
perfect Hessenberg matrices integer conjugate to each other, see
Example~\ref{2reduced}. This happens also for matrices in
$SL(2,\z)$.

\subsection{Perfect Hessenberg matrices conjugate to a given
one}\label{toHess} In this subsection we show the algorithm to
construct a perfect Hessenberg matrices for a given
$SL(n,\z)$-matrix. It is based on the following proposition.

\begin{proposition}\label{DavHes}
Let $M$ be an $SL(n,\z)$-matrix with irreducible over $\q$
characteristic polynomial. For any integer primitive vector $v$
$($i.e., with relatively prime coordinates$)$ there exists a
unique matrix $C$ such that

--- $C(e_1)=v$;

--- the matrix $CMC^{-1}$ is perfect Hessenberg $($we denote this matrix by $(M|v)$$)$.
\end{proposition}

\begin{remark}
This means that for any lattice preserving linear operator any
primitive integer vector can be extended to the basis of integer
lattice in a way such that the matrix of the operator is perfect
Hessenberg in this basis.
\end{remark}

\begin{proof}
First, we construct the corresponding perfect Hessenberg matrix.
Let $M$ be a matrix in $SL(n,\z)$ with irreducible characteristic
polynomial, and $A$ be a linear operator with matrix $M$ in some
fixed integer basis. Take any primitive integer vector $v$  and
consider a set of vector spaces
$$
V_i= \sspan\big(v, A(v), A^2(v),\ldots, A^{i-1}(v)\big),
$$
for $i=1,\ldots,n$ (here we denote the span of vectors $v_1,
\ldots ,v_m$ by $\sspan(v_1, \ldots , v_m)$). Since the
characteristic polynomial of $A$ is irreducible, the dimension of
$V_i$ equals $i$ and the set of all spaces $V_i$ forms a complete
flag in $\r^n$. Since for any integer $j$ the vector $A^j(v)$ is
integer, the spaces $V_i$ contain an integer sublattice of rank
$i$.

\vspace{1mm}

Let us inductively construct an integer basis $\{e_i\}$ of a
vector space $\r^n$ such that:

---  for $i=1,\ldots, n$, the vectors $e_1,\ldots e_i$ form a basis
of the integer sublattice $\z^n\cap V_i$;

--- the matrix of the operator $A$ is perfect Hessenberg in this basis.

\vspace{1mm}

{\parindent 0pt {\it Base of induction}. We put $e_1=v$. It is
clear that the vector $e_1$ generate the one-dimensional integer
sublattice of $V_1$.}

\vspace{1mm}

{\parindent 0pt {\it Step of induction.} Suppose we have
constructed $e_i$ for all $i\le k$, such that $e_1,\ldots, e_k$ is
a basis of integer sublattice contained in $V_i$. Let us find
$e_{k+1}$.}

By construction, $e_1,\ldots, e_k$ generate the integer sublattice
$\z^n\cap V_k$. Hence there exists an integer vector $g_{k+1}$
such that $e_1,\ldots, e_k,g_{k+1}$ generate the integer
sublattice $\z^{n}\cap V_{k+1}$. Actually $g_{k+1}$ is one of the
integer primitive vectors of $V_{k+1}$ with the smallest possible
nonzero distance to the space $V_k$ (all such vectors are
contained in two hyperplanes parallel to $V_k$).

Since $A(e_k)$ is contained in $V_{k+1}$, it is decomposable in
the basis $e_1,\ldots, e_k, g_{k+1}$ with integer coefficients:
$$
A(e_k)=\sum\limits_{i=1}^{k}q_{i,k} e_i+a_{k+1,k}g_{k+1}.
$$
For $i=1,\ldots, k$ we define $b_{i,k}$, and $a_{i,k}$ as integer
quotients and reminders:
$$
q_{i,k}=b_{i,k}\cdot|a_{k+1,k}|+a_{i,k},
$$
where $0\le a_{i,k} < a_{k+1,k}$. Rewrite
$$
A(e_k)=|a_{k+1,k}|
\left(\sign(a_{k+1,k})g_{k+1}+\sum\limits_{i=1}^{k}b_{i,k}
e_i\right)+\sum\limits_{i=1}^{k}a_{i,k} e_i.
$$
satisfying $0\le a_{i,k} < a_{k+1,k}$ for $i=1,\ldots,k$. Finally
we put
$$
e_{k+1}=\sign(a_{k+1,k})g_{k+1}+\sum\limits_{i=1}^{k}b_{i,k} e_i.
$$
Since the characteristic polynomial of $A$ is irreducible, the
integer spaces $V_i$ are not invariant subspaces of $A$, hence
$e_1, \ldots, e_{k+1}$ are linearly independent and generate
$V_{k+1}$.

The matrix $\hat M$ of the operator $A$ in the basis $\{e_i\}$ is
of Hessenberg type
$$
\Big\langle a_{1,1}, |a_{2,1}|\Big|a_{1,2}, a_{2,2},
|a_{3,2}|\Big|
\cdots
\Big|a_{1,n-1},\ldots, a_{n-1,n-1}, |a_{n,n-1}|\Big\rangle.
$$
By the definition, $\hat M$ is a perfect Hessenberg matrix. The
matrices $M$ and $\hat M$ represent the same operator $A$ in two
different integer bases, hence $M$ and $\hat M$ are integer
conjugate.

Since the characteristic polynomial of $A$ is irreducible, the
integer spaces $V_i$ are not invariant subspaces of $A$. Hence,
the integers $a_{i+1,i}$ are nonzero for $i=1,\ldots, n{-}1$.
Therefore, the Hessenberg complexity of $\hat M$ is positive.

Denote by $C$ the transition matrix to the basis $\{e_i\}$. Then
we have  $C(e_1)=v$ and $\hat M=CMC^{-1}$ is perfect Hessenberg.

\vspace{2mm}

Finally, we say a few words about uniqueness of $\hat M$. The
spaces $V_i$ are uniquely defined. The vector $e_1$ is uniquely
defined. On each step there is a unique way to define $e_{k+1}$.
Hence the transition matrix $C$ is uniquely defined. Therefore,
such matrix $\hat M$ is unique.
\end{proof}

Let us briefly outline the algorithm used in the proof of
Proposition~\ref{DavHes}.
\\
{\bf Algorithm to construct perfect Hessenberg matrices}.

 {\it Input Data.} We are given by a matrix $M$ of a lattice preserving
operator $A$ with irreducible characteristic polynomial and an
integer vector $v$.

{\it Step 1.} We put $e_1=v$.

{\it Inductive Step k.} Suppose we have constructed $e_i$ for all
$i\le k$. For $g_{k+1}$ we take one of the integer primitive
vectors of $V_{k+1}$ with the smallest possible nonzero distance
to the space $V_k$. Find the coordinates $q_{i,k}$ for
$i=1,\ldots, k$ and $a_{k+1,k}$ from the decomposition
$$
A(e_k)=\sum\limits_{i=1}^{k}q_{i,k} e_i+a_{k+1,k}g_{k+1}.
$$
For $i=1,\ldots, k$ find $b_{i,k}$, and $a_{i,k}$ as integer
quotients and reminders:
$$
q_{i,k}=b_{i,k}\cdot|a_{k+1,k}|+a_{i,k}.
$$
Then we have
$$
e_{k+1}=\sign(a_{k+1,k})g_{k+1}+\sum\limits_{i=1}^{k}b_{i,k} e_i.
$$
Finally, let $C$ be a transition matrix to the basis $\{e_k\}$.

{\it Output Data.}  In the output we have the perfect Hessenberg
matrix $CMC^{-1}$.

\vspace{2mm}

We use the following corollary in the proof of Theorem~\ref{cfr} below.

\begin{corollary}\label{orbitinvariance}
Consider an $SL(n,\z)$-operators $A$ with matrix $M$. Let $B$ be
an arbitrary $GL(n,\z)$-operator commuting with $A$. Then for an
arbitrary $v$ we have
$$
(M|v)=(M|B(v)).
$$
\end{corollary}

\begin{proof}
Since $A$ commutes with $B$, each step of the above algorithm is
invariant produces the same data for both $v$ and $B(v)$. Hence
$(M|v)=(M|B(v))$.
\end{proof}

\subsection{Reduction to $\varsigma$-reduced Hessenberg matrices}\label{Hess1}

In this subsection we show the existence of $\varsigma$-reduced
Hessenberg matrices integer conjugate to a given one. We show how
to find it explicitly in Subsection~\ref{reconstruct} after we
introduce Klein-Voronoi continued fractions.

\begin{theorem}\label{toRedHess}
For any matrix $M$ in $SL(n,\z)$ with irreducible characteristic
polynomial  there exists a $\varsigma$-reduced Hessenberg matrix
$\tilde M$ with positive Hessenberg complexity such that $M$ is
integer conjugate to $\tilde M$.
\end{theorem}

\begin{proof} By Proposition~\ref{toHess} there exists at least one perfect
Hessenberg matrix integer conjugate to $M$. Since the set of
values of Hessenberg complexity is discrete and bounded from
below, there exists a perfect Hessenberg matrix $\tilde M$ integer
conjugate to~$M$ and with minimal possible Hessenberg complexity.
By definition  $\tilde M$ is a $\varsigma$-reduced Hessenberg
matrix.
\end{proof}

\subsection{Finiteness of $\varsigma$-reduced Hessenberg matrices}\label{Hess2}
In this subsection we prove the following theorem.

\begin{theorem}\label{finiteness}
For any $SL(n,\z)$-matrix $M$ with irreducible characteristic
polynomial there exists finitely many $\varsigma$-reduced
Hessenberg matrices integer conjugate to $M$.
\end{theorem}

In the proof of this theorem we use the following general
proposition.

\begin{proposition}\label{type&char}
Any Hessenberg matrix with positive Hessenberg complexity is
uniquely defined by its Hessenberg type and the characteristic
polynomial. \qed
\end{proposition}

\begin{proof}
Consider a Hessenberg matrix $M=(a_{i,j})$ of a given Hessenberg
type with positive Hessenberg complexity. From the Hessenberg type
of $M$ we know all its columns except for the last one. Let us
show that the last column in
 is uniquely defined by the
coefficients of its characteristic polynomial. Let the
characteristic polynomial of $M$ be
$$
x^n+c_{n-1}x^{n-1}+\cdots+c_1x+c_0.
$$
Direct calculations show that for any $k$ the coefficient $c_k$ is
a polynomial in $a_{i,j}$ variables that does not depend on
$a_{1,n},\ldots, a_{k,n}$. The unique monomial for $c_{k}$
containing $a_{k+1,n}$ is
$$
\left(\prod\limits_{j=k+1}^{n-1}a_{j{+}1,j}\right)a_{k+1,n}.
$$
Since the Hessenberg complexity of $M$ is nonzero, the product in
the brackets is nonzero. Hence $a_{k+1,n}$ is a function of $c_k$
and the elements $a_{i,j}$ contained in the first $n-1$ columns.
This concludes the proof of the proposition.
\end{proof}

The following example shows that simply Hessenberg complexity
together with characteristic polynomial do not distinguish all the
integer conjugacy classes.

\begin{example}
The following two matrices
$$
\left(
\begin{array}{ccc}
0 &1 &3\\
1 &0 &0\\
0 &3 &8\\
\end{array}
\right) \quad \hbox{and} \quad
\left(\begin{array}{ccc}
0 &2 &5\\
1 &1 &2\\
0 &3 &7\\
\end{array}
\right)
$$
are not integer conjugate but have the same Hessenberg complexity
equal to 3 and the same characteristic polynomials.
\end{example}

\vspace{2mm}

{\it Proof of Theorem~\ref{finiteness}.} The existence of
$\varsigma$-reduced Hessenberg matrices integer conjugate to $M$
follows from Theorem~\ref{toRedHess}. By definition they all have
the same Hessenberg complexity (say, $c$). The number of
Hessenberg types whose Hessenberg complexity equals $c$ is finite.
It is clear that the integer conjugate matrices have the same
characteristic polynomial, hence by Proposition~\ref{type&char}
there exists at most one Hessenberg matrix of a given Hessenberg
type integer conjugate to $M$. Therefore, there is only a finite
number of $\varsigma$-reduced Hessenberg matrices integer
conjugate to $M$.
\qed

\subsection{Families of Hessenberg matrices with given
Hessenberg type}\label{Struc}

Denote by $H(\Omega)$ the set of all Hessenberg matrices in
$SL(n,\z)$ of Hessenberg type $\Omega$.

For an arbitrary Hessenberg type
$$
\Omega=\langle
a_{1,1},a_{1,2}|a_{2,1},a_{2,2},a_{2,3}|\cdots|a_{n{-}1,1},\ldots,
a_{n{-}1,2},a_{n{-}1,n}\rangle
$$
and $k=1,\ldots, n{-}1$ we denote by $v_k(\Omega)$ the vector
$(a_{k,1},\ldots,a_{k,k+1},0,\ldots,0)$, and by $M_k(\Omega)$ ---
the matrix with zero first $n{-}1$ columns and the last one equals
to $v_k(\Omega)$.

Denote by $\sigma(\Omega)$ the $(n{-}1)$-dimensional simplex with
vertices $O$, $O{+}v_1, \ldots, O{+}v_{n-1}$ where $O$ is the
origin.

\begin{definition}
The {\it integer volume} of a simplex $\sigma$ with integer
vertices is the index of the sublattice generated by the edges of
$\sigma$ in the lattice of all integer vectors in the plane
spanned by $\sigma$.
\end{definition}

\begin{theorem}\label{family}
Let $\Omega$ be a Hessenberg type.

{\it i$)$.} The set $H(\Omega)$ is not empty if and only if the
integer volume of $\sigma(\Omega)$ equals one.

{\it ii$)$.} Suppose that $M_0\in H(\Omega)$, then $H(\Omega)$ is
an integer affine $(n{-}1)$-dimensional sublattice in the lattice
of all integer $(n\times n)$-matrices. More precisely,
$$
H(\Omega)=\left\{ M_{0}+\sum\limits_{i=1}^{n-1} c_i M_i(\Omega)
\Big| c_1, \ldots, c_{n-1} \in \z \right\}.
$$
\end{theorem}

The proof of Theorem~\ref{family} is based on
Lemma~\ref{structLemma} which we show after the following
definition.

\begin{definition}
Consider an integer vector $v$ and a $k$-dimensional plane $\pi$
containing the integer sublattice of rank $k$ such that $v$ is not
in $\pi$. The {\it integer distance} from $v$ to $\pi$ is the
index of the sublattice generated by the integer vectors of the
set $\pi\cup\{u\}$ in the whole integer lattice of the
$(k{+}1)$-dimensional plane spanning $v$ and $\pi$.
\end{definition}

\begin{lemma}\label{structLemma}
Consider a Hessenberg matrix $M$ of type $\Omega$, let its last
column be an integer vector vector $v$. The matrix $M$ is in
$SL(n,\z)$ if and only if the following conditions hold:

--- the integer volume of $\sigma(\Omega)$ equals one;

--- the integer distance from the vector $v$ to the integer
hyperplane containing $\sigma(\Omega)$ equals one.
\end{lemma}

\begin{proof}
Let $A$ be an operator with Hessenberg matrix $M$ in the basis $\{
e_i\}$ of integer lattice.

Suppose that $M$ is in $SL(n,\z)$, then the operator $A$ preserves
all integer volumes and integer distances. Since the integer
volume of the coordinate $(n{-}1)$-dimensional simplex $S_e^{n-1}$
with vertices
$$
O, O{+}e_1, \ldots, O{+}e_{n-1}
$$
equals one, the integer volume of the image
$\sigma(\Omega)=A(S_e^{n-1})$ equals one. Notice that
$$
A(\sspan(S_e^{n-1}))=\sspan(\sigma(\Omega)) \qquad \hbox{and}
\qquad A(e_n)=v.
$$
Since the integer distance from the point $O{+}e_n$ to the plane
spanned by the vectors $e_1, \ldots, e_{n-1}$ equals one, the
integer distance from the point $O+v$ to the integer hyperplane
containing $\sigma(\Omega)$ also equals one.

Suppose now that both conditions of the lemma hold. Then the
operator $A$ takes the integer lattice (generated by $ e_1,\ldots,
e_n$) to itself bijectively. Therefore, $M$ is in $SL(n,\z)$.
\end{proof}

{\it Proof of Theorem~\ref{family}.} {\it $($i$)$}  Suppose the
integer volume of $\sigma(\Omega)$ equals one. Then we choose $v$
to be at unite integer distance to the plane
$\sspan(\sigma(\Omega))$. Then by Lemma~\ref{structLemma} we get
the matrix. Conversely if $H(\Omega)$ contains an
$SL(n,\z)$-matrix, then by Lemma~\ref{structLemma} the integer
volume of $\sigma(\Omega)$ equals 1.

\vspace{1mm}

Statement~{\it$($ii$)$} is straightforward, since the determinant
of the matrix is additive with respect to the operation of
addition of vectors in the last column. \qed

 \vspace{2mm}

We conclude this subsection with a particular example.

\begin{example}\label{ex1}
Let us consider matrices of Hessenberg type $\langle 0,1|1,0,2
\rangle$. All matrices of that type form a two-parametric family
$$
H(\langle 0,1|1,0,2 \rangle)=\left\{
\left(
\begin{array}{ccc}
0 &1 &1\\
1 &0 &0\\
0 &2 &1\\
\end{array}
\right) +
m
\left(
\begin{array}{ccc}
0 &0 &0\\
0 &0 &1\\
0 &0 &0\\
\end{array}
\right) +
n
\left(
\begin{array}{ccc}
0 &0 &1\\
0 &0 &0\\
0 &0 &2\\
\end{array}
\right)\Bigg| m,n \in \z \right\}.
$$
We denote
$$
H_{\langle 0,1|1,0,2\rangle}^{(1,0,1)}(m,n)= \left(
\begin{array}{ccc}
0 &1 &n+1\\
1 &0 &m\\
0 &2 &2n+1\\
\end{array}
\right)
$$
with integer parameters $m$ and $n$. The discriminant of the
matrix $H_{\langle 0,1|1,0,2\rangle }^{(1,0,1)}(m,n)$ equals
$$
-44-44n^2-56mn-32n^3+32m^3+16m^2n^2+16mn^2+16m^2n-56n-8m+52m^2.
$$

The set of matrices with negative discriminant for the given
family coincides with the union of integer solutions of the
following inequalities:
$$
2m\le -n^2-n-2 \quad \hbox{and} \quad 2n\le m^2+m.
$$
\begin{figure}
$$\epsfbox{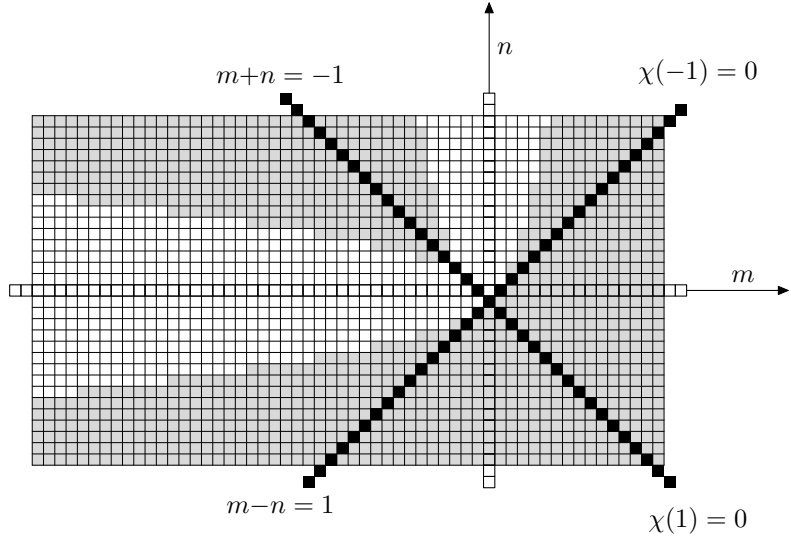}$$
\caption{The family of matrices of Hessenberg type $\langle
0,1|1,0,2\rangle$.}\label{result.3}
\end{figure}

In Figure~\ref{result.3} the square in the intersection of the
$m$-th column and the $n$-th row corresponds to the matrix
$H_{\langle 0,1|1,0,2 \rangle}^{(1,0,1)}(m,n)$. Black squares
correspond to the matrices with reducible characteristic
polynomials. Light gray squares correspond to the matrices with
three real eigenvalues. The rest have a pair of complex conjugate
eigenvalues.
\end{example}

\section{Complete geometric invariant of conjugacy classes}\label{KlVor}

In this section we introduce a geometric complete invariant of
integer conjugacy classes: multidimensional continued fractions in
the sense of Klein-Voronoi. We start with the classical
two-dimensional case corresponding to Gauss Reduction Theory in
Subsection~\ref{GRT}, where we show a relation between integer
conjugacy classes of matrices and certain continued geometric
continued fractions. We give general definitions of Klein-Voronoi
continued fractions in all dimensions in Subsection~\ref{cf_1}.
Klein-Voronoi continued fractions for matrices of $SL(n,\z)$
possess additional combinatorial periodicity, we discuss it in
Subsection~\ref{cf_3}. Finally in Subsection~\ref{cf_4} we show
that Klein-Voronoi continued fractions classify integer conjugacy
classes of $SL(n,\z)$-matrices.

\subsection{Geometry of Gauss Reduction Theory}\label{GRT}

In this subsection we briefly discuss geometry of two-dimensional
case. We skip all the proofs in this subsection, for a more
detailed exposition with all proofs we refer
to~\cite{Karpenkov2010} and~\cite{Karpenkov2008}.

\vspace{2mm}

It is usual to split $SL(2,\z)$ in the following three cases.

\vspace{2mm}

{\bf Complex case}: a characteristic polynomial of such matrices
has two complex conjugate roots. There are only three classes of
such matrices, they are represented by
$$
\left(
\begin{array}{cc}
1&1\\
-1&0\\
\end{array}
\right), \quad
 \left(
\begin{array}{cc}
0&1\\
-1&0\\
\end{array}
\right), \quad
\hbox{and} \quad
\left(
\begin{array}{cc}
0&1\\
-1&-1\\
\end{array}
\right).
$$

{\bf Degenerate case}: the characteristic polynomial has a double
root (that is actually equal to 1). Such matrices are integer
conjugate to exactly one of the following family
$$
\left(
\begin{array}{cc}
1&n\\
0&1\\
\end{array}
\right)\quad \hbox{for $n\ge 0$.}
$$

{\bf Totally real case}: the case of two real roots is the most
complicated. The geometric description is given via continued
fractions. So we give necessary definitions at first.

\subsubsection{Ordinary continued fractions} The
expression (finite or infinite)
$$
a_0+1/(a_1+1/(a_2+\ldots)\ldots))
$$
is an {\it ordinary continued fraction} if $a_0 \in \z$, $a_k\in
\z_+$ for $k> 0$. Denote it by $[a_0: a_1; \ldots]$ (or by $[a_0:
a_1; \ldots ;a_n]$).

An ordinary continued fraction is {\it odd $($even$)$} if it has
an odd (even) number of elements. Any rational number has a unique
odd and even ordinary continued fractions:
$$
\frac{7}{5}=1+\frac{1}{2+\frac{\displaystyle 1}{\displaystyle
2}}=1+\frac{1}{2+\frac{\displaystyle 1}{\displaystyle 1+1/1}}.
$$
The odd and even continued fractions of the same number coincide
except for the very last elements, as in example with 7/5:
$$
[a_0: a_1; \ldots ;a_n]=[a_0: a_1; \ldots ;a_n{-}1:1].
$$
Any irrational number has a unique infinite ordinary continued
fraction.

\subsubsection{Integer geometry notation}
Let us briefly recall some notions of integer geometry. A point is
{\it integer} if all its coordinates are integers. A segment or a
vector is integer if it has integer endpoints. An angle is {\it
integer} if its vertex is integer and its edges contain integer
points distinct to the vertex.

\begin{definition}
The {\it integer length} of an integer segment $AB$ is the number
of inner integer points in the segment plus one, we denote it by
$\il(AB)$.
\\
The {\it integer sine} of an integer angle $ABC$ is the index of
the sublattice generated by integer vectors of the edges of the
angle in the lattice of integer points, we denote it by $\isin
(ABC)$.
\end{definition}

For additional information on lattice trigonometry we refer
to~\cite{Karpenkov2008} and~\cite{Karpenkov2009}.

\begin{definition}
Consider an arbitrary angle $C$ with vertex at the origin. The
boundary of the convex hull of all integer points in $C$ except
for the origin is called the {\it sail} for $C$.
\end{definition}
In general a sail is a broken line that contains a finite or
infinite number of vertices.

\begin{definition}
Consider an arbitrary angle with integer vertex. Let the sail for
this angle be a broken line with the sequence of vertices $(V_i)$.
Denote:
$$
\begin{array}{l}
a_{2k}=\il V_kV_{k+1},\\
a_{2k-1}=\isin V_{k-1}V_{k}V_{k+1}\\
\end{array}
$$
for all admissible indices. The {\it lattice length-sine sequence}
({\it LLS-sequence}, for short) for the sail is the sequence
$(a_n)$.
\end{definition}

\subsubsection{Geometry of ordinary continued fractions}
An odd or infinite continued fraction of any real number
$\alpha\ge 1$ has the following geometric interpretation. Consider
the angle in the first orthant defined by two rays $y=\alpha x$
and $y=0$, we denote it by $C_\alpha$. Let also the first vertex
of the sail for $C_\alpha$ be in the ray $y=0$ (actually it is the
point $(1,0)$).

\begin{theorem}
Consider a real number $\alpha\ge 1$. Let $(a_0,\ldots, a_{2n})$
$($or $(a_0,a_1,\ldots)$$)$ be the LLS-sequence of $C_\alpha$.
Then
$$
\alpha=[a_0: a_1; \ldots; a_{2n}] \qquad
(\alpha=[a_0:a_1;\ldots]).
$$
\qed
\end{theorem}

We refer to~\cite{Karpenkov2011a} for the geometry of continued
fractions in a more general situation.

\begin{figure}
{\large
$$
\begin{array}{c}
\epsfbox{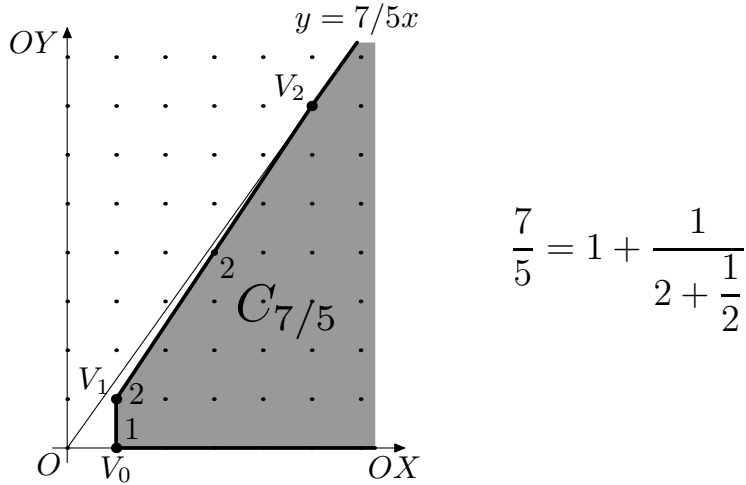}
\end{array} \qquad
\frac{\displaystyle 7}{\displaystyle 5}=1+\frac{\displaystyle
1}{\displaystyle 2+\frac{\displaystyle 1}{\displaystyle 2}}
$$
} \caption{Geometry of the continued fraction for
$7/5=[1:2;2]$.}\label{7_5.1}
\end{figure}

\begin{example}
On Figure~\ref{7_5.1} we show the example of the sail for
$C_{7/5}$. The boundary convex hull consists of two rays and two
segments. It contains three vertices $V_0$, $V_1$, and~$V_2$.
Direct calculations show that
$$
\il(V_0V_1)=1; \quad \isin(V_0V_1V_2)=2; \quad \hbox{and} \quad
\il(V_1V_2)=2.
$$
So we have $7/5=[1:2;2]$.
\end{example}

\vspace{2mm}

\subsubsection{Continued fractions for 2-dimensional operators with
real eigenvectors}  Let us expand a notion of continued fractions
in the following way. Consider an arbitrary $GL(2,\r)$-operator
with two real distinct eigenvalues. This operator has exactly two
eigenlines. The complement to these eigenlines consists of four
angles. The boundary of the convex hull of all integer points
except the origin inside any of these angles is called a {\it
sail} of the matrix. The set of all four sails is called the {\it
geometric continued fraction}.

\begin{figure}
$$\epsfbox{defcf.1}$$
\caption{The periodic sails of a matrix
${{7\hbox{ } 18}\choose {5 \hbox{ } 13}}$
%
%
are periodic. One of its periods is $(1,1,3,2)$.}\label{matrix}
\end{figure}

Let $A$ be an $SL(2,\z)$-operator with two real distinct
eigenvector. Then all four sails of the operators are two-side
infinite broken lines. Moreover all four LLS-sequences coincide up
to a shift and/or reversing the order. An operator $A$ acts on its
sails as a shift, therefore, all the LLS-sequences of $A$ are
periodic and the period is defined by the shift.

\begin{example}
On Figure~\ref{matrix} we show an example of a geometric continued
fraction for the matrix
$$
\left(
\begin{array}{cc}
7&18\\
5&13\\
\end{array}
\right).
$$
The LLS-sequences of all four sails are periodic. Their periods
are either $(1,1,3,2)$ or $(2,3,1,1)$ counterclockwise.
\end{example}

\subsubsection{Totally real case} Note that, there are two subfamilies of totally real
$SL(2,\z)$-operators: with positive eigenvectors and with negative
eigenvalues. The composition with an operator of symmetry about
the origin (i.e.~$-Id$) gives a one-to-one correspondence between
the sets of all operators with positive and negative eigenvalues.
Hence we restrict ourselves to the case of operators with positive
eigenvalues.

\begin{definition}
A matrix
$\left(
\begin{array}{cc}
a&c\\
b&d\\
\end{array}
\right)$
in $SL(2,\z)$ is {\it reduced} if $d> b \ge a \ge 0$.
\end{definition}

It is interesting that the LLS-sequence of a reduced operator can
be written directly from the coefficients of the matrix.

\begin{theorem}\label{ggg}
Consider a reduced matrix $M={{a\hbox{ } c}\choose {b \hbox{ }
d}}$. Suppose $\frac{b}{a}=[a_1;a_2:\ldots :a_{2n-1}]$ and
$\lambda=\lfloor \frac{d-1}{b} \rfloor$ then
$$
(a_1,a_2,\ldots ,a_{2n-1},\lambda)
$$
is one of the periods of the LLS-sequences for the geometric
continued fraction of $M$. \qed
\end{theorem}

We say that the reduced matrices are ''almost'' normal forms,
since each matrix could have more than one normal form. Their
number is described via LLS-sequence is as follows.

\begin{theorem}{\bf(On almost normal forms.)}
The number of reduced matrices in an integer conjugacy class with
minimal period $(a_1,\ldots, a_k)$ of the corresponding
LLS-sequence is $k$. \qed
\end{theorem}

{\noindent {\it Remark.} Let us say a few words about a relation
between reduced and $\varsigma$-reduced matrices in $SL(2,\z)$.
From one hand any reduced totally real matrix is a perfect
Hessenberg matrix. From the other hand any $\varsigma$-reduced
matrix is also reduced, although there are certain reduced
matrices that are not $\varsigma$-reduced.}

\vspace{2mm}

The LLS-sequence itself is the complete invariant of the set of
all integer conjugacy classes in $SL(2,\z)$.

\begin{theorem}\label{complete2D}{\bf(On complete invariant of integer conjugacy classes.)}
An even period of the LLS-sequence $($up to an even shift and
reversing the order$)$ is a complete invariant of an integer
conjugacy class of a $SL(2,\z)$-matrices with distinct positive
real eigenvalues. \qed
\end{theorem}

Notice that all odd periods correspond to matrices with negative
determinant.

\begin{example}
For the matrix
$$M=\left(
\begin{array}{cc}
1519&1164\\
-1964&-1505\\
\end{array}
\right)
$$
the period is $(1,2,1,2)$. So there are two reduced matrices: with
the period $(1,2,1,2)$ and $(2,1,2,1)$. The coefficients $a$ and
$b$ for the matrices are then respectively as follows
$$
\frac{b}{a}=[1;2:1]=\frac{4}{3} \quad \hbox{or} \quad
\frac{b}{a}=[2;1:2]=\frac{8}{3}.
$$
Now it is not hard to find the elements $c$ and $d$ of the reduced
matrices from conditions $\lambda=\lfloor \frac{d-1}{b} \rfloor$
and $ad-bc=1$. Finally we get all two reduced matrices integer
conjugate to $M$. They are as follows:
$$
\left(
\begin{array}{cc}
3&8\\
4&11\\
\end{array}
\right)
\quad \hbox{and} \quad
\left(
\begin{array}{cc}
3&4\\
8&11\\
\end{array}
\right).
$$
The complexities of these matrices are 4 and 8 respectively. So
the first matrix is $\varsigma$-reduced.
\end{example}

\subsubsection{Matrices with the same Hessenberg type, their geometric continued fractions}

It turns out that matrices of the same Hessenberg type has similar
geometric background. We start with the following example.

\begin{example}
Consider the family of matrices of Hessenberg type ${\langle
3,4\rangle}$:
$$
\left(
\begin{array}{cc}
3&3m+2\\
4&4m+3\\
\end{array}
\right).
$$
For $m=-1,-2$ the corresponding operators does not have real
eigenvectors. For $m\ge 0$ the periods of the corresponding
geometric continued fractions are
$$
[2,2], \quad [1,2,1,1], \quad [1,2,1,2], \quad [1,2,1,3], \quad
\ldots, \quad [1,2,1,m], \quad \ldots
$$
The periods in cases $m=0$ and $m=2$ are twice the minimal. This
means that the corresponding matrices are squares of some other
integer matrices. For $m=0$ it is a square of some integer matrix
with negative determinant, since the minimal period is odd, and
for $m=2$ it is a square of some $SL(2,\z)$-matrix with negative
determinant. The matrices are $\varsigma$-reduced for $m\ge 2$.

For $m\le -3$ the periods of the corresponding geometric continued
fractions are
$$
[4,1], \quad [4,2], \quad [4,3], \quad [4,4], \quad \ldots, \quad
[4,-2-m], \quad \ldots
$$
For $m=-6$ the matrix is a square of some $GL(2,\z)$-matrix.
Starting from $m\le -6$ the matrices are $\varsigma$-reduced.
\end{example}

So almost all matrices of Hessenberg type ${\langle 3,4\rangle}$
are $\varsigma$-reduced. This is actually the case for all
Hessenberg types in $SL(2,\z)$.
\begin{theorem}
Almost all matrices of a given Hessenberg type in $SL(2,\z)$ are
$\varsigma$-reduced. \qed
\end{theorem}
This follows from general Theorem~\ref{cfr} below and direct
calculation of complexities for all vertices of the period, we
skip the proof here.

For further information about the two-dimensional case we refer
the reader, for instance, to the works~\cite{Karpenkov2010},
\cite{Lewis1997} and~\cite{Manin2002}.

\subsection{Continued fractions in the sense of Klein-Voronoi}

\subsubsection{Background}

In 1839 C.~Hermite~\cite{Hermite1839} posed the problem of
generalizing ordinary continued fractions to the
higher-dimensional case. Since then there were many different
definitions generalizing different properties of ordinary
continued fractions. A nice geometrical generalization of ordinary
continued fraction for operators with all real eigenvalues was
made by F.~Klein in~\cite{Klein1895} and~\cite{Klein1896}.

Multidimensional continued fractions in the sense of Klein have
many relations with other branches of mathematics. For example,
O.~N.~German~\cite{German2002} and
J.-O.~Moussafir~\cite{Moussafir2000} discussed the connection
between the sails of multidimensional continued fractions and
Hilbert bases. In~\cite{Tsuchihashi1983} H.~Tsuchihashi described
the relationship between periodic multidimensional continued
fractions and multidimensional cusp singularities.
M.~L.~Kontsevich and Yu.~M.~Suhov studied the statistical
properties of random multidimensional continued fractions
in~\cite{Kontsevich1999}. The relations to approximation theory of
maximal commutative subgroups is discussed by A.~Vershik and the
author in~\cite{Karpenkov2010b}. The combinatorial topological
generalization of Lagrange theorem was obtained by E.~I.~Korkina
in~\cite{Korkina1994} and its algebraic generalization by
G.~Lachaud~\cite{Lachaud1993}. The book~\cite{Arnold2002} of
V.~I.~Arnold is a good survey of geometric problems and theorems
associated with one-dimensional and multidimensional continued
fractions in the sense of Klein (see also his
articles~\cite{Arnold1989}, \cite{Arnold1998},
and~\cite{Arnold1999}).

Approximately at the time of the works by F.~Klein G.~Voronoi in
his dissertation~\cite{Voronoui1952} introduced a geometric
algorithmic definition for all the cases even for operators with
pairs of complex conjugate eigenvalues. In~\cite{Buchmann1985}
and~\cite{Buchmann1985a} J.~A.~Buchmann generalized Voronoi's
algorithm making it more convenient for computation of fundamental
units in orders. We use ideas of J.~A.~Buchmann to define the {\it
multidimensional continued fraction in the sense of Klein-Voronoi}
for all the cases. Note that if all the eigenvalues of an operator
are real numbers then the Klein-Voronoi multidimensional continued
fraction is a continued fractions in the sense of Klein.

\subsubsection{General definitions}\label{cf_1}

Consider an operator $A$ in $GL(n,\r)$ with distinct eigenvalues.
Suppose that it has $k$ real eigenvalues $r_1,\ldots, r_k$ and
$2l$ complex conjugate eigenvalues $c_1,\bar c_1,\ldots, c_l, \bar
c_l$, where $k+2l=n$.

Denote by $T^l(A)$ the set of all real operators commuting with
$A$ such that their real eigenvalues are all unit and the absolute
values for all complex eigenvalues equal one. Actually, $T^l(A)$
is an abelian group with operation of matrix multiplication.

For a vector $v$ in $\r^n$ we denote by $T_A(v)$ the orbit of $v$
with respect of the action of the group of operators $T^l(A)$. If
$v$ is in general position with respect to the operator $A$ (i.e.
it does not lie in invariant planes of $A$), then $T_A(v)$ is
homeomorphic to the $l$-dimensional torus. For a vector of an
invariant plane of $A$ the orbit $T_A(v)$ is also homeomorphic to
a torus of positive dimension not greater than $l$, or to a point.

\begin{example}
Suppose that $A$ is a totally real operator. Since all its
eigenvectors are real,  $T^0(A)$ consists only of the unit
operator and $T_A(v)=\{v\}$.
\end{example}

\begin{example}
Now consider an operator $A$ with a pair of complex eigenvalues
whose all the other eigenvalues are real. The group $T^1(A)$
correspond to elliptic rotations in the invariant plane of $A$
corresponding to complex eigenvalues. Such rotations are
parameterized by an angle of rotation. A general orbit of $T_A(v)$
is an ellipse around the $(n{-}2)$-dimensional invariant subspace
corresponding to real eigenvalues. Any orbit in the invariant
subspace of real eigenvalues consists of one point.
\end{example}

Let $g_i$ be a real eigenvector with eigenvalue $r_i$ for
$i=1,\ldots, k$; $g_{k+2j-1}$ and $g_{k+2j}$ be vectors
corresponding to the real and imaginary parts of some complex
eigenvector with eigenvalue $c_j$ for $j=1,\ldots, l$. Consider
the coordinate system corresponding to the basis $\{g_i\}$:
$$
OX_1X_2\ldots X_{k}Y_{1}Z_{1}Y_2Z_2\ldots Y_l Z_l.
$$

Denote by $\pi$ the $(k{+}l)$-dimensional plane $OX_1X_2\ldots X_k
Y_1Y_2\ldots Y_l$. Let $\pi_+$ be the cone in the plane $\pi$
defined by the equations $y_i\ge 0$ for $i=1,\ldots, l$. For any
$v$ the orbit $T_A(v)$ intersects the cone $\pi_+$ in a unique
point.

\begin{definition}
A point $p$ in the cone $\pi_+$ is said to be {\it $\pi$-integer}
if the orbit $T_A(p)$ contains at least one integer point.
\end{definition}

Consider all (real) hyperplanes invariant under the action of the
operator $A$. There are exactly $k$ such hyperplanes. In the above
coordinates the $i$-th of them is defined by the equation $x_i=0$.
The complement to the union of all invariant hyperplanes in the
cone $\pi_+$ consists of $2^k$ arcwise connected components.
Consider one of them.

\begin{definition}
The convex hull of all $\pi$-integer points except the origin
contained in the given arcwise connected component is called a
{\it factor-sail} of the operator $A$. The set of all factor-sails
is said to be the {\it factor-continued fraction} for the operator
$A$.
\\
The union of all orbits $T_A(*)$ in $\r^n$ represented by the
points in the factor-sail is called the {\it sail} of the operator
$A$. The set of all sails is said to be the  {\it continued
fraction} for the operator $A$ in the sense of Klein-Voronoi (see
in Figure~\ref{complex.1} below).
\end{definition}

The intersection of the factor-sail with a hyperplane in $\pi$ is
said to be an $m$-dimensional face of the factor-sail if it is
homeomorphic to the $m$-dimensional disc.

The union of all orbits in $\r^n$ represented by points in some
face of the factor-sail is called the {\it orbit-face} of the
operator $A$.

Integer points of the sail are said to be {\it vertices} of this
sail.

\subsubsection{Algebraic continued fractions}\label{cf_3}

Consider now an operator $A$ in the group $GL(n,\z)$ with
irreducible characteristic polynomial. Suppose that it has $k$
real roots $r_1,\ldots, r_k$ and $2l$ complex conjugate roots:
$c_1,\bar c_1,\ldots, c_l, \bar c_l$, where $k+2l=n$. In the
simplest possible cases $k{+}l=1$ any factor-sail of $A$ is a
point. If $k{+}l>1$, than any factor-sail of $A$ is an infinite
polyhedral surface homeomorphic to $\r^{k+l-1}$.

\begin{definition}
The group of all $GL(n,\z)$-operators commuting with $A$ is called
the {\it Dirichlet group} and denoted by $\Xi(A)$.\\
The subgroup of the Dirichlet group $\Xi(A)$ consisting of all
matrices whose real eigenvalues are all positive is called the
{\it positive Dirichlet group}. We denote it by $\Xi_+(A)$.
\end{definition}


The Dirichlet group $\Xi(A)$ takes the Klein-Voronoi continued
fraction to itself but maybe exchange the sails. The positive
Dirichlet group $\Xi_+(A)$ consists exactly from operators
preserving all the sails. By Dirichlet unit theorem (see, for
instance, in~\cite{Borevich1966}) the group $\Xi (A)$ is
homomorphic to $\z^{k+l-1}\oplus G$, where $G$ is some finite
commutative group. The group $\Xi_+(A)$ is homeomorphic to
$\z^{k+l-1}$ and its action on any sail is free. The quotient of a
sail by the action of $\Xi_+(A)$ is homeomorphic to the
$(n{-}1)$-dimensional torus.

\begin{definition}\label{defFD}
A {\it fundamental domain of the Klein-Voronoi continued fraction}
is a collection of open orbit-faces such that for any
$\Xi(A)$-orbit of orbit-faces of the continued fraction there
exists a unique representative in the collection.
\\
A {\it fundamental domain of a sail} is a collection of open
orbit-faces such that for any $\Xi_+(A)$-orbit of orbit-faces of
the sail there exists a unique representative in the collection.
\end{definition}

\begin{example}
Let us study an operator $A$ with a Frobenius matrix
$$
\left(
\begin{array}{ccc}
0 &0 &1\\
1 &0 &1\\
0 &1 &3\\
\end{array}
\right).
$$
This operator has one real and two complex conjugate eigenvalues.
Therefore, the cone $\pi_+$ for $A$ is a two-dimensional
half-plane. In Figure~\ref{complex.1}a the halfplane $\pi_+$ is
colored in light gray and the invariant plane corresponding to the
pair of complex eigenvectors is in dark gray. The vector shown in
Figure~\ref{complex.1}a with endpoint at the origin is an
eigenvector of~$A$.

In Figure~\ref{complex.1}b we show the cone $\pi_+$. The invariant
plane separates $\pi_+$ onto two parts. The dots on $\pi_+$ are
the $\pi$-integer points. The boundaries of the convex hulls in
each part of $\pi_+$ are two factor-sails. Actually, one
factor-sail is taken to another by the induced action of $-Id$,
where $Id$ is an identity operator of $\r^3$.

Finally, in Figure~\ref{complex.1}c we show one of the sails.
Three orbit-vertices shown in the figure correspond to the vectors
$(1,0,0)$, $(0,1,0)$, and $(0,0,1)$: the large dark points
$(0,1,0)$ and $(0,0,1)$ are visible on the corresponding
orbit-vertices.

The positive Dirichlet group $\Xi_+(A)$ in our example is
homeomorphic to $\z$, it is generated by $A$. The group $\Xi(A)$
is homeomorphic to $\z\oplus\z/2\z$ with generators $A$ and $-Id$.
The operator $A$ takes the point $(1,0,0)$ and its orbit-vertex to
the point $(0,1,0)$ and the corresponding orbit-vertex. Therefore,
a fundamental domain of the continued fraction for the operator
$A$ contains one orbit-vertex and one vertex edge. For instance,
we can choose the orbit-vertex corresponding to the point
$(1,0,0)$ and the orbit-edge corresponding to the ''tube''
connecting orbit-vectors for the points $(1,0,0)$ and $(0,1,0)$.

\begin{figure}
$$\epsfbox{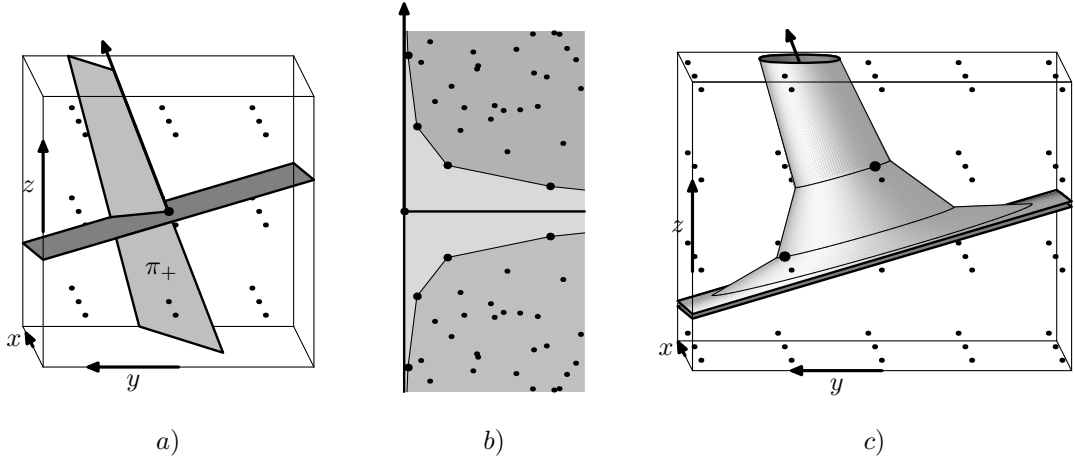}$$
\caption{A tree-dimensional example: a) the cone $\pi_+$ and the
eigenplane; b) the continued factor-fraction; c) a sail of the
continued fraction.}\label{complex.1}
\end{figure}
\end{example}

\subsection{Invariants of conjugacy classes}\label{cf_4}
As we already know, any integer conjugacy class of totally-real
$SL(2,\z)$-matrices is defined by a period of the corresponding
geometric continued fraction (Theorem~\ref{complete2D}). Similar
statement is true in the multidimensional case.

\subsubsection{Theorem on complete invariant}
Let $A$ be a $GL(n,\z)$-operator with distinct eigenvalues. Then
any $B\in \Xi(A)$ acts on its Klein-Voronoi continued fraction of
$A$ as a periodic shift.

\begin{definition}
Let $P$ be a Klein-Voronoi continued fraction (of some
$GL(n,\z)$-opera/-tor). A transformation $T$ of $P$ is said to be
a {\it period} of $P$ if there exists a $GL(n,\z)$-operator $A$
such that

--- the Klein-Voronoi continued fraction of $A$ coincides with $P$;

--- the transformation $T$ coincides with the action of $A$ on $P$.

\end{definition}

Let us define now the congruence for continued fractions and their
periods.

\begin{definition}
Two Klein-Voronoi continued fractions are said to be {\it integer}
congruent if there exists a $GL(n,\z)$-operator that takes one of
them to another.
\\
Two periods of the congruent Klein-Voronoi continued fractions are
said to be {\it integer congruent} if there exists a
$GL(n,\z)$-operator that takes the first sail to the second and
the period of the first sail to the period of the second.
\end{definition}

Periods of Klein continued fractions in a totally-real three
dimensional case $(k=3)$ were studied in works~\cite{Korkina1996},
\cite{Korkina1995}, \cite{Lachaud1993}, \cite{Lachaud2002},
\cite{Karpenkov2004}, \cite{Karpenkov2004a}, etc.

\begin{theorem}\label{teorem_inv}{\bf(On complete invariant in general case.)}
A Klein-Voronoi continued fraction together with one of its
periods is a complete invariant of an integer  conjugacy class of
$GL(n,\z)$-matrices.
\end{theorem}

\begin{remark}
In Theorem~\ref{PeriodStructure} below we give the description of
the set of all periods of a Klein-Voronoi continued fraction.
\end{remark}

\begin{proof}
If two matrices are integer conjugate then their Klein-Voronoi
continued fractions and the corresponding periods are integer
congruent.

Suppose now that two matrices $A$ and $B$ have integer congruent
continued fractions with congruent periods (let $C$ be the integer
congruence of the sails and periods). Denote
$$
\tilde B=CBC^{-1}.
$$
Matrices $A$ and $\tilde B$ have the same periods of the same
continued fraction. Therefore, their actions coincide for all
points of the Klein-Voronoi sail for $A$. It is clear that points
of the Klein-Voronoi continued fractions span $\r^n$, hence by
linearity $A=\tilde B$.
\end{proof}

\subsubsection{Structure of the set of periods of Klein-Voronoi continued fraction}
Consider a Klein-Voronoi continued fractions $S$. The composition
of two integer linear shifts of $S$ is again an integer shift of
$P$ and hence the set of all periods is the group. Denote it by
$\Xi(P)$.

\begin{theorem}\label{PeriodStructure}
Let $A$ be a matrix with irreducible characteristic polynomial
over $\q$ and $P$ be its Klein-Voronoi continued fraction. Then
the group of its periods $\Xi(P)$ coincides with the Dirichlet
group $\Xi(A)$.
\end{theorem}

Here we should show that distinct Dirichlet groups does not
correspond to the same Klein-Voronoi continued fraction. So the
proof of the theorem is based on the following proposition.

\begin{proposition}\label{invariant}
Consider operators $A$ and $B$ in $SL(n,\z)$ with irreducible
characteristic polynomials. Operators $A$ and $B$ commute if and
only if they have the same Klein-Voronoi continued fraction.
\end{proposition}

{\noindent {\it Remark.} Recall that Klein-Voronoi continued
fractions are defined only for operators with distinct
eigenvalues. If the characteristic polynomial of an operator is
irreducible over $\q$, then all its roots are distinct, so any
operator of Proposition~\ref{invariant} has a Klein-Voronoi
continued fraction.}

\begin{proof}
If operators $A$ and $B$ with irreducible characteristic
polynomials commute, then they have the same eigenvectors. Hence
an arbitrary operator $C$ commutes with $A$ if and only if $C$
commutes with $B$. Hence $T_A(v)=T_B(v)$ for any vector $v$.
Therefore, the Klein-Voronoi continued fractions of both operators
coincide by construction.

\vspace{2mm}

Let us prove the converse statement. Let the Klein-Voronoi
continued fractions for $A$ and $B$ coincide as sets.

Suppose $A$ has real eigenvectors $e_1,\ldots, e_k$ and complex
conjugate eigenvectors $a_j\pm Ib_j$ for $j=1,\ldots,l$, where
$k+2l=n$ (here $I=\sqrt{-1}$). Let $g_i$ be a real eigenvector
with eigenvalue $r_i$ for $i=1,\ldots, k$; $g_{k+2j-1}$ and
$g_{k+2j}$ be vectors corresponding to the real and imaginary
parts of some complex eigenvector with eigenvalue $c_j$ for
$j=1,\ldots, l$. Consider the coordinate system corresponding to
the basis $\{g_i\}$:
$$
OX_1X_2\ldots X_{k}Y_{1}Z_{1}Y_2Z_2\ldots Y_l Z_l.
$$
In this coordinates we consider the form $\Phi_A$ that in the
above coordinates is written as
$$
\Phi_A(x_1,\ldots, x_n)=\left( \prod\limits_{i=1}^k
x_i\prod\limits_{j=1}^l(y_j^2+z_j^2) \right)
$$
(we study this form later in Section~\ref{algorithm}). Similarly
we define the form $\Phi_B$ for the operator $B$. From definition
it follows that $A$ preserves $\Phi_A$ and $B$ preserves $\Phi_B$.

Since asymptotically (in the complement to the balls centered at
the origin with the radius increasing to infinity) the
Klein-Voronoi continued fraction for $A$ (for $B$) is tends to the
set $\Phi_A=0$ (and $\Phi_B=0$ respectively) in a continuous
category, the operators $A$ and $B$ have all the same invariant
subspaces. In particular, their one-dimensional real eigenspaces
corresponding to real eigenvectors and two-dimensional eigenspaces
(we denote them by $\pi_1$, \ldots, $\pi_l$) defined by pairs of
complex conjugate roots coincide. This implies that $A$ and $B$
commute if and only if they commute for the vectors of the
invariant planes $\pi_1$, \ldots, $\pi_l$.

\vspace{1mm}

Let us show that $A$ and $B$ restricted to the plane $\pi_i$
($i=1,\ldots, l$) commute. Consider the section of the
Klein-Voronoi continued fraction by the plane passing through $v$
and parallel to the invariant subspace spanned by all complex
eigenvectors of $A$, we denote it by $T_v$. From construction
$T_v=T_A(v)$. By the above the invariant subspace spanned by all
complex eigenvectors of $A$ coincide with invariant subspace
spanned by all complex eigenvectors of $B$, and hence
$T(v)=T_B(v)$. Therefore, $T(v)=T_A(v)=T_B(v)$.

It is clear the forms $\Phi_A$ and $\Phi_B$ are constant at the
orbit $T(v)$ for any $v$. This follows from the fact that any
operator $P$ of $T^l(A)$ (respectively $T^l(B)$) preserves
$\Phi_A$ (respectively, $\Phi_B$), since all eigenvalues of $P$
are of unit absolute values and $P$ is diagonalizable in the
eigenbasis of $A$ (respectively $B$). Therefore, from linearity
reasoning $\Phi_A=c\cdot\Phi_B$ for some nonzero constant $c$.
Therefore, the operator $A$ preserves the form $\Phi_B$.

Consider now the plane $\pi_j$ for some $1\le j \le l$ and take
coordinates $OXY$ such that the restriction of the form
$\Phi_B(v)$ to this plane is
$$
x^2+y^2=\lambda.
$$
Direct calculations show that there are two types of operators
that preserve this form. They are written in $OXY$ coordinates as
follows
$$
\left(
\begin{array}{cc}
\alpha &\beta\\
-\beta & \alpha\\
\end{array}
\right) \quad \hbox{and} \quad
\left(
\begin{array}{cc}
-\alpha &\beta\\
\beta & \alpha\\
\end{array}
\right)
$$
with parameters $\alpha$ and $\beta$. The operators of the second
family have two real eigenvalues in the plane $\pi_j$, which is by
the above not the case for the operator $B$. Therefore, both $A$
and $B$ are from the first family. All operators of the first
family commute. Hence $A$ commutes with $B$ in planes $\pi_j$ for
$1\le j \le l$.

Therefore, $A$ and $B$ are both diagonalisable in the same complex
basis. Hence $A$ and $B$ commute.
\end{proof}

{\noindent {\it Remark.} It is interesting to notice that for
certain $A$ there exist some operators corresponding to the second
family in the proof of Proposition~\ref{invariant}. These
operators preserve the Klein-Voronoi continued fraction of $A$,
although they have different Klein-Voronoi continued fractions or
even multiple roots.}

\vspace{2mm}

{\noindent{\it Proof of Theorem~\ref{PeriodStructure}.} From one
hand, by Proposition~\ref{invariant} all $SL(n,\z)$-matrices
defining the shift commute with each other and also with $A$.
Therefore, they are in $\Xi(A)$. From the other hand, any matrix
of $\Xi(A)$ defines a period of Klein-Voronoi continued fraction
$P$. \qed}

\section{Algorithmic aspects of reduction to $\varsigma$-reduced
matrices}\label{algorithm}

In Section~\ref{HMCC} above we show the existence of and
finiteness of $\varsigma$-reduced matrices in each integer
conjugacy class of $SL(n,\z)$-matrices. The aim of this section is
to show the techniques to construct $\varsigma$-reduced matrices
integer conjugate to a given one. In Subsection~\ref{cf_2} we give
a geometric interpretation of the Hessenberg complexity as a
volume of a certain simplex, which is called the
MD-characteristics. Further in Subsection~\ref{cf_5} we use the
MD-characteristics to show that that all $\varsigma$-reduced
matrices are obtained from integer vertices of Klein-Voronoi
continued fraction by applying the algorithm of
Subsection~\ref{toHess} to them. The corresponding techniques is
discussed in Subsection~\ref{reconstruct}.

\subsection{Markoff-Davenport characteristics}\label{cf_2}
In this Subsection we characterize the Hessenberg complexity in
terms of Markoff-Davenport characteristics.

\subsubsection{Definition of the MD-characteristics and its invariance under the action of the
Dirichlet group} The study of the Markoff-Davenport
characteristics is closely related to the theory of minima of
absolute values of homogeneous forms with integer coefficients in
$n$-variables of degree~$n$. One of the first works in this area
was written by A.~Mar\-koff~\cite{Markoff1879} for the
decomposable forms (into the product of real linear forms) for
$n=2$. Further, H.~Davenport in series of
works~\cite{Davenport1938},~\cite{Davenport1938a},~\cite{Davenport1939},
\cite{Davenport1941}, and~\cite{Davenport1943} made first steps
for the case of decomposable forms for $n=3$.

\vspace{2mm}

Consider $A\in SL(n,\z)$. Denote by $P(A,v)$ the parallelepiped
spanned by vectors $v$, $A(v)$, $\ldots$, $A^{n-1}(v)$, i.e.,
$$
P(A,v)=\bigg\{O+\sum\limits_{i=0}^{n-1}\lambda_iA^{i}(v)\bigg|0\le
\lambda_i \le 1, i=0,\ldots, n{-}1\bigg\},
$$
where $O$ is the origin.
\begin{definition}
The {\it Markoff-Davenport characteristics} (or {\it
MD-characteristic}, for short) of an $SL(n,\z)$-operator $A$ is a
functional:
$$
\Delta_A:\r^n\to \r \qquad \hbox{defined by}
\qquad\Delta_A(v)=V(P(A,v)),
$$
where $V(P(A,v))$ is the nonoriented volume of $P(A,v)$.
\end{definition}

\begin{proposition}
Consider $A\in SL(n,\z)$ and let $B\in \Xi(A)$. Then for an
arbitrary $v$ we have
$$
\Delta_A(v)=\Delta_A(B(v)).
$$
\end{proposition}

{\noindent {\it Remark.} This means that the MD-characteristics
naturally defines a function over the set of all orbits of the
Dirichlet group.}

\begin{proof}
Since $B\in\Xi(A)$ we have $A^nB(v)=BA^n(v)$. Hence the
parallelepiped $P(A,B(v))$ coincides with $B(P(A,v))$. Since $B\in
SL(n,\z)$, the volume of the parallelepiped is preserved.
Therefore,
$$
\Delta_A(v)=\Delta_A(B(v)).
$$
\end{proof}

\subsubsection{Homogeneous forms associated to $SL(n,
\z)$-operators}

Let $\{e_i\}$ be an integer basis of $\r^n$. Consider any
$SL(n,\z)$-operator $A$ with irreducible characteristic
polynomial. Suppose that it has $k$ real eigenvalues $r_1,\ldots,
r_k$ and $2l$ complex conjugate eigenvalues $c_1,\bar c_1,\ldots ,
c_l, \bar c_l$, where $k+2l=n$. Let us now define a new basis of
vectors $g_1, \ldots, g_{k+2l}$ in the following way. For
$i=1,\ldots, k$ we choose $g_i$ to be an eigenvector corresponding
to the eigenvalue $r_i$. For $j=1,\ldots, l$ we choose
$g_{k+2j-1}$ and $g_{k+2j}$ to be the real and the imaginary parts
of some complex eigenvector corresponding to the eigenvalue $c_j$.
Consider the system of coordinates
$$
OX_1X_2\ldots X_{k}Y_{1}Z_{1}Y_2Z_2\ldots Y_l Z_l
$$
corresponding to the basis $\{g_i\}$.

A form
$$
\alpha \left( \prod\limits_{i=1}^k
x_i\prod\limits_{j=1}^l(y_j^2+z_j^2) \right)
$$
with nonzero $\alpha$ is said to be {\it associated} to the
operator $A$.

\vspace{2mm}

\begin{proposition}\label{DavChar}
Let $A$ be an $SL(n,\z)$-operator with irreducible characteristic
polynomial. Then the MD-characteristics of $A$ is an absolute
value of a form associated to $A$ for a certain nonzero $\alpha$.
\end{proposition}

\begin{proof}
Let us consider the formulas of MD-characteristics of $A$ in the
eigen-basis of vectors
$$
g_1, \ldots, g_k, g_{k+1}{+}Ig_{k+2}, g_{k+1}{-}Ig_{k+2}, \ldots,
g_{k+2l-1}{+}Ig_{k+2l}, g_{k+2l-1}{-}Ig_{k+2l}
$$
in $\c^n$, where $I=\sqrt{-1}$. Let the coordinates in this
eigen-basis be $\{t_i\}$.

Then for any vector $v=(t_1,\ldots,t_n)$ we have
$$A^j(x)=(r_1^j t_1,\ldots, r_k^j t_k,c_1^j t_{k+1}, \bar c_1^j t_{k+2},
\ldots, c_l^j t_{k+2l-1}, \bar c_l^j t_{k+2l}).
$$
Therefore,
$$
\Delta_A(t_1,\ldots,t_n)=
\alpha\left|\prod\limits_{i=1}^k t_i\prod\limits_{j=1}^l
(t_{k+2j-1}t_{t+2j})\right|=
\frac{\alpha}{4^l} \left| \prod\limits_{i=1}^k
x_i\prod\limits_{j=1}^l(y_j^2+z_j^2) \right|
$$

Simple calculations show that $\alpha\ne 0$.
\end{proof}

\subsubsection{Hessenberg complexity in terms of MD-characteristics}\label{DavDir}

\begin{proposition}\label{DavHes0}
Consider an operator $A$ with Hessenberg matrix $M$ in some
integer basis $\{e_i\}$. The Hessenberg complexity $\varsigma(M)$
equals the value of MD-characteristics $\Delta_A(e_1)$.
\end{proposition}

\begin{proof}
Suppose that the Hessenberg type of the matrix $M$ is
$$
\big\langle a_{1,1},a_{1,2}\big|a_{2,1},a_{2,2},a_{2,3}\big|\cdots
\big|a_{n{-}1,1},\ldots, a_{n{-}1,2},a_{n{-}1,n}\big\rangle.
$$
Denote by $V_{k}$ the plane $\sspan \big(v, A(v), A^2(v),\ldots,
A^{k-1}(v) \big)$.

Let us inductively show that
$$
A^k(e_1)=\left(\prod_{i=1}^ka_{i,i+1}\right)e_{k+1} +v_k \quad
\hbox{where $v_k\in V_{k}$}.
$$

{\noindent {\it Base of induction.} We have
$A(e_1)=a_{1,2}e_2+a_{1,1}e_1$.}

\vspace{1mm}

{\noindent {\it Step of induction.} Suppose that the statement
holds for $k=m$, i.e.,
$$
A^m(e_1)=\left(\prod_{i=1}^ma_{i,i+1}\right)e_{m+1} +v_m, \quad
\hbox{and $v_m\in V_{m}$}.
$$
Let us show the statement for $m{+}1$. Since $M$ is Hessenberg,
$A(v_m)$ is in $V_{m+1}$. Therefore, we have
$$
\begin{array}{lcl}
A^{m+1}(e_1)&=&
A\left(\left(\prod\limits_{i=1}^ma_{i,i+1}\right)e_{m+1}\right)+A(v_m)=
\\
&&\left(\prod\limits_{i=1}^{m+1}a_{i,i+1}\right)e_{m+1}+
\left(A(v_m) +\left(\prod\limits_{i=1}^ma_{i,i+1}\right) \big(
A(e_{m+1}){-}a_{m+1,m+2}e_{m+2}\big)\right).
\end{array}
$$

The second summand in the last expression is in $V_{m+1}$. We have
shown the step of induction.}

Therefore,
$$
\Delta_A(e_1)=\prod\limits_{i=1}^{n{-}1}|a_{i{+}1,i}|^{n-i}=\varsigma(M).
$$
This concludes the proof of the proposition.
\end{proof}
Further we use the following corollary.
\begin{corollary}\label{toHessCor}
Consider an operator $A$ with Hessenberg matrix $M$ in some
integer basis $\{e_i\}$. Let $v$ be any primitive integer vector.
Then we have
$$
\varsigma(M|v)=\Delta_A(v),
$$
where $(M|v)$ is the matrix constructed by the algorithm of
Subsection~\ref{toHess}. \qed
\end{corollary}

\subsection{Klein-Voronoi continued fractions and minima of MD-characteristics}\label{cf_5}

In the following theorem we use Klein-Voronoi continued fractions
to find minima of MD-characteristics.

\begin{theorem}\label{cfr}
Consider a matrix $M\in SL(n,\z)$ with distinct eigenvalues. Let
$U$ be a fundamental domain of the Klein-Voronoi continued
fractions for $M$ $($see Definition~\ref{defFD}$)$. Then we have:

{\it $($i$)$} For any $\varsigma$-reduced matrix $\hat M$ integer
conjugate to $M$ there exists $v\in U$ such that $ \hat M=(M|v)$.

{\it $($ii$)$} Let $v\in U$. The matrix $(M|v)$ is
$\varsigma$-reduced if and only if the MD-characteristic
$\Delta_A(v)$ attains its minimal value.
\end{theorem}

\begin{proof}
Notice that Theorem~\ref{cfr}$($ii$)$ is a direct corollary of
Corollary~\ref{toHessCor}.

\vspace{2mm}

Let us prove Theorem~\ref{cfr}$($i$)$. Let $A$ be an
$SL(n,\z)$-operator with irreducible characteristic polynomial
defined by the matrix $M$. By Proposition~\ref{DavChar} there
exists a nonzero constant $\alpha$ such that MD-characteristics
$\Phi_A$ at any point in the system of coordinates $OX_1X_2\ldots
X_{k}Y_{1}Z_{1}Y_2Z_2\ldots Y_l Z_l$ is
$$
\alpha \left| \prod\limits_{i=1}^k
x_i\prod\limits_{i=1}^l(y_i^2+z_i^2) \right|.
$$
Suppose that the minimal absolute value of $F$ on the set of
integer points except the origin equals $m_0$.

Choose the coordinates $OX_1\ldots X_kY_1Y_2,\ldots,Y_l$ in the
cone $\pi_+$. Consider a projection of $\r^n$ to the cone along
the $T_A(v)$ orbits. Since we project along the $T_A(v)$ orbits on
which the MD-characteristics is constant, the projection of the
MD-characteristics is well-defined, denote it by $\tilde\Phi_A$.
In the chosen coordinates of $\pi_+$, the function $\tilde\Phi_A$ is
written as follows:
$$
\alpha\big|\prod_{j=1}^kx_j\prod_{j=1}^ly_j^2 \big|.
$$
The obtained function is convex in any orthant of the cone
$\pi_+$. Since any factor sails are boundary of a convex hull in each orthant,
all the minima of the convex function $\tilde\Phi_A$ restricted to the convex hulls
are attained at the boundary, i.e., at $\pi$-integer points of factor-sails.
Therefore, all integer minima of $\Phi_A$ are at vertices of the Klein-Voronoi
continued fraction.

By Corollary~\ref{toHessCor}, the Hessenberg complexity
$\varsigma(M|v)$ coincides with MD-characteristics $\Delta_A(v)$.
Since any matrix integer conjugate to $M$ has a presentation in
the form $(M|v)$ and all the integer minima of MD-characteristics
attained at vertices of the Klein-Voronoi continued fraction, any
$\varsigma$-reduced operator $\tilde M$ is represented as
$(M|v_0)$ for some vertex $v_0$ of Klein-Voronoi continued
fraction. By Corollary~\ref{orbitinvariance}, for any $B\in\Xi(A)$
we have:
$$
(M|B(v_0))=(M|v_0),
$$
since all such $B$ commutes with $A$.
Hence a vector $v_0$ can be chosen from the fundamental domain
$D$. This concludes the proof.
\end{proof}

Let us give an example of two $\varsigma$-reduced perfect
Hessenberg matrices integer conjugate to each other.

\begin{example}\label{2reduced}
The $\varsigma$-reduced Hessenberg matrices (with Hessenberg
complexity equal to 3)
$$
M_1=\left(
\begin{array}{ccc}
0 &1 &2\\
1 &0 &0\\
0 &3 &5\\
\end{array}
\right) \qquad \hbox{and} \quad
M_2=\left(\begin{array}{ccc}
0 &2 &3\\
1 &1 &1\\
0 &3 &4\\
\end{array}
\right)
$$
are integer conjugate.

The reason for this is as follows. Consider the Klein-Voronoi
continued fraction of $A$ with matrix $M_1$. It contains integer
vertices $p_1=(1,0,0)$ and $p_2=(0,1,-1)$. It turns out that $p_1$
and $p_2$ are not in the same orbit of the Dirichlet group but
have the same MD-characteristics equals~3. Hence we get distinct
two integer conjugate $\varsigma$-reduced Hessenberg matrices:
$M_1=(M_1|(1,0,0))$ and $M_2=(M_1|(0,1,-1))$.
\end{example}

\subsection{Construction of $\varsigma$-reduced
matrices by Klein-Voronoi continued fraction}\label{reconstruct}
Any $\varsigma$-reduced Hessenberg matrix for the operator $A$ is
constructed starting from some vertex in a fundamental domain of
the Klein-Voronoi multidimensional continued fraction as follows.

\vspace{2mm} {\bf Techniques to find $\varsigma$-reduced matrices
in an integer conjugacy class.}

{\it Step 1.} Find a fundamental domain of the Klein-Voronoi
continued fraction for the operator $A$ (see
Remark~\ref{construction}).

{\it Step 2.} Take all vertices of the constructed fundamental
domain and find among them all vertices with minimal value of the
MD-characteristics (say $v_1,\ldots, v_k$).

{\it Step 3.} By Theorem~\ref{cfr}$($i$)$ and~$($ii$)$ all the
$\varsigma$-reduced matrices integer conjugate to $M$ are
$(M|v_1),\ldots, (M|v_k)$. They are all constructed by the
algorithm described in Subsection~\ref{toHess}.

\vspace{2mm}

\begin{remark}
For the case of $SL(2,\z)$ we have geometric Gauss Reduction
Theory: each vertex of geometric continued fraction corresponds to
a starting point of the period. The corresponding matrix is
written according to Theorem~\ref{ggg}.
\end{remark}

\begin{remark}\label{construction}
Currently Step~1 is the most complicated. The totally real case of
matrices with all eigenvalues being real is studied quite good.
For the algorithms of constructing multidimensional continued
fractions in this case, we refer to the papers by
R.~Okazaki~\cite{Okazaki1993},
J.-O.~Moussafir~\cite{Moussafir2000a} and the
author~\cite{Karpenkov2009b}. E.~Korkina in~\cite{Korkina1996}
and~\cite{Korkina1995}, G.~Lachaud
in~\cite{Lachaud1993},~\cite{Lachaud2002}, A.~D.~Bruno and
V.~I.~Parusnikov in~\cite{Bryuno1994}, ~\cite{Parusnikov1999},
~\cite{Parusnikov2000} and ~\cite{Parusnikov2005} the author
in~\cite{Karpenkov2004} and~\cite{Karpenkov2004a} produced a large
number of fundamental domains for periodic algebraic
two-dimensional continued fractions (see also the
site~\cite{Briggs} by K.~Briggs). Some fundamental domains in
three dimensional case are found for instance
in~\cite{Karpenkov2006a}. The case with complex conjugate
eigenvalues is relatively new, we are planning to study it in our
forthcoming paper.
\end{remark}

\begin{example}
Let us consider an example of an operator $A$ defined by the
matrix
$$
\left(
\begin{array}{ccc}
-2&-4 &-3\\
1&2 & 2\\
-1&-1 & 3\\
\end{array}
\right).
$$
The characteristic polynomial of this operator has three distinct
real roots. Therefore, the Klein-Voronoi continued fraction
consists of 8 sails. The compositions of operators $-Id$, $A$, and
$2Id+A^{-1}$ define equivalence between all these sails (here $Id$
is the identity operator) and hence all $\varsigma$-reduced operators can
be written from vertices of one sail. Consider a sail
containing the point $[1,0,0]$. There are exactly three distinct orbits of
the Dirichlet group containing the vertices in this sail. They are
defined by the following points
$$
[0,0,1], \quad [1,0,0],\quad \hbox{and} \quad [3,-1,1].
$$
(We skip all the calculations of convex hulls corresponding to the
sail, see the algorithms in~\cite{Moussafir2000a},
\cite{Karpenkov2009b}, \cite{Shintani1976}, \cite{Okazaki1993}).
The MD-characteristics of these vectors are respectively: $1$,
$2$, and $4$. So the minimum of the MD-characteristics (which is
$1$ in this case) is attained on the vertices of the orbit of the
Dirichlet group containing $[0,0,1]$ (and on the corresponding
orbits for the rest seven sails). Therefore, there exists a unique
$\varsigma$-reduced Hessenberg matrix, which is
$$
\left(
\begin{array}{ccc}
0&0 & 1\\
1&0 & 1\\
0&1 &-3\\
\end{array}
\right).
$$
The perfect Hessenberg matrices for the vertices $[1,0,0]$ and
$[3,-1,1]$ are respectively
$$
\left(
\begin{array}{ccc}
0&1 & -1\\
1&0 & 0\\
0&2 &-3\\
\end{array}
\right) \quad \hbox{and} \quad \left(
\begin{array}{ccc}
1&0 &-1\\
2&0 & 3\\
0&1 &-4\\
\end{array}
\right),
$$
their $\varsigma$-complexities are $2$ and $4$.
\end{example}

\bibliographystyle{plain}      
\bibliography{sl3zz}

\vspace{5mm}

\end{document}